\documentclass[12pt]{amsart}
\usepackage{epic,eepic}
\usepackage[dvips]{graphics,color}
\newtheorem{theorem}{Theorem}[section]
\newtheorem{proposition}[theorem]{Proposition}
\newtheorem{lemma}[theorem]{Lemma}

\newtheorem{remark}[theorem]{Remark}

\newtheorem{conjecture}[theorem]{Conjecture}

\newcommand{\BB}[1]{\mathbb{#1}}
\newcommand{\CAL}[1]{\mathcal{#1}}

\newcommand{\BF}[1]{\mathbf{#1}}
\newcommand{\BFIT}[1]{\mbox{\boldmath$\mathit{#1}$}}

\renewcommand{\descriptionlabel}[1]%
{\hspace\labelsep \upshape\bfseries #1}

\newcommand{\SCBFIT}[1]{\mbox{\boldmath\scriptsize$\mathit{#1}$}}
\newcommand{\COMP}{\raisebox{0.3ex}{\hspace{0.3ex}%
{$\scriptstyle{\circ}$}\hspace{0.5ex}}}   %%% composition
\newcommand{\Hdeg}{\operatorname{Hdeg}}
\newcommand{\Det}{\operatorname{Det\,}}
\newcommand{\Span}{\operatorname{Span\,}}

\newcommand{\Bul}{\circle*{.5}}
\newcommand{\QED}{\hspace*{\fill}\raisebox{.54ex}%
{$\framebox[.35em]{\/}$}\\ \hspace*{4ex}}
\title{Restrictions of Smooth Functions\\
to a Closed Subset}%
\author{Shuzo IZUMI}
\date{\today}
\begin{document}
%%%%%%%%%%%%%%%%%%%%%%%%%%%%%%%%%%%%%%%
\setcounter{page}{1}
\maketitle
%%%%%%%%%%%%%%%%%%%%%%%%%%%%%%%%%%%%%%%%%%%%%%%%%%
We first provide an approach to the recent conjecture of 
Bierstone-Milman-Paw\l ucki on Whitney's old problem 
on $C^d$ extendability of functions defined on a closed subset 
of a Euclidean space, using the higher order paratangent bundle 
they introduced. 
For example, the conjecture is affirmative for classical 
fractal sets. 
Next, we give a sharpened form of Spallek's theorem 
on controllability of flatness by the values on a closed set. 
The multi-dimensional Vandermonde matrix plays an important role 
in both cases. 
%%%%%%%%%%%%%%%%%%%%%%%%%%%%%%%%%%%%%%%%%%%%%%%%%%
\section*{Introduction}
%%%%%%%%%%%%%%%%%%%%%%%%%%%%%%%%%%%%%%%%%%%%%%%%%%%

There remain very simple unsolved problems in a 
close neighbourhood of elementary calculus. 
They are concerned with the values 
of smooth functions on a subset as follows. 
Let $X$ be a closed subset of a domain $\Omega$ of a Euclidean 
space $\BB{R}^n$. 
\begin{enumerate}
\item
\textbf{Extension problem}:
\textit{
Find the condition for a function on $X$ to be extendable 
to a $C^d$ function on $\Omega$.
}
\item
\textbf{Flatness problem}:
\textit{
If $X$ is thick enough around $\BF{\xi}$, the Taylor expansion 
of a differentiable function $f$ defined in a neighbourhood of 
$\BF{\xi}$ is determined up to some order. 
Find a geometric expression of this thickness.
}
\end{enumerate}

Whitney \cite{whitney-2} posed the extension problem 
and gave a necessary and sufficient condition in the case 
$n=1$. 
Glaeser \cite{g1} solved this problem in the case $d=1$ 
introducing the linearized paratangent bundle. 
(Originally, the spelling is ``paratingent" 
in French. cf. \cite{bouligand} and \cite{g1}.) 
Recently, Bierstone, Milman and Paw\l ucki \cite{bmp2} 
introduced a very interesting geometric notion 
``higher order paratangent bundle'' $\tau_N^d(X)$ 
to this problem, generalizing Glaeser's paratangent bundle. 
This is a variant of higher order tangent bundle of $X$. 
Further, they associated to every 
function $f:X\longrightarrow \mathbb{R}$, a subbundle 
$\nabla_N^d f \subset \tau_N^d(X)\times \mathbb{R}$ 
$(d,\ N \in \mathbb{N})$ over $X$. 
In these constructions of $\tau_N^d(X)$ and $\nabla_N^d f$, 
they used a set operation found by Glaeser (see \S \ref{ptg}). 
They posed a general conjecture on the extension problem 
as follows. \\
%%%%%%%%%%%%%%%%%%%%%%%
$(*)$ \textit{
A function $f$ is $C^d$ extendable if and only if $\nabla_N^d f$ 
is the graph of a mapping of $\tau_N^d(X)$ into $\mathbb{R}$ 
for a suitable $N$.
}\\
%%%%%%%%%%%%%%%%%%%%

They proved that, for the closure $X$ of 
an open subset of a regular submanifold of a Euclidean domain, 
$\tau_N^d(X)$ coincides with the full higher order tangent 
bundle of the submanifold on $X$ and that the conjecture $(*)$ is 
affirmative for such $X$. 
Further, they have obtained a positive result for the important 
case of compact subanalytic sets. 

We call the second problem flatness problem 
since uniqueness up to order $r$ of Taylor expansion at 
$\xi \in X$ is assured by controlling $r$-flatness 
by the values on $X$. 
This problem was first considered by Spallek \cite{spallek}. 
The present author \cite{jam} developed 
fundamental properties of ``Spallek function", an invariant 
defined for the germ $X_{\xi}$ which measures the efficiency of 
flatness control of functions at $\xi$ by their values on $X$. 

Since Whitney's works, it has been widely known that 
the theory of smooth functions is closely related 
to the interpolation theory (e.g. \cite{whitney-2}, 
\cite{interpolation}, \cite{kergin}, \cite{mm}, \cite{spallek}, 
\cite{bmp2}, \cite{jam}). The reason is that differential properties 
are not punctual but ``molecular" (Glaeser) as seen in the 
bi-punctual inequality used to define Whitney function 
(\cite{whitney-1}). 
Glaeser \cite{interpolation} proposed two methods 
of application of interpolation, Lagrange interpolation and 
``interpolation schemes" to treat differential properties. 
He put emphasis on the latter. But 
we adopt the former method in this paper. 
Following Glaeser, we treat 
interpolations with $(n+d)!/n!d!$ nodes for the problems of 
$C^d$ functions on $\mathbb{R}^n$. The most important point 
is that ``Vandermonde matrix" appears in the matrix 
representation of Taylor expansion 
(see the proofs of \ref{ext} and \ref{flat}). 
We measure the decrease of the Vandermonde determinants of 
accumulating nodal sets. 

In \ref{ext} we show the following. 
Suppose that $X \subset \mathbb{R}^n$ includes 
nodal sets accumulating to a point and that their Vandermonde 
determinants are not rapidly decreasing relative to their diameters. 
Then $X$ has the full higher order paratangent spaces at 
the accumulating point. 
Remember that if $X$ has full higher order paratangent spaces 
at any point of $X$, the conjecture $(*)$ is valid for such $X$ 
(see \cite{bmp2}). Many classical fractal sets, such as 
Cantor set, Koch curve, Sierpinski gasket and Menger sponge, 
satisfy this condition (see \ref{self-similar}). 

As for flatness, we need a more quantitative argument. 
We give a sufficient condition \ref{flat} for sets $X$ 
to control flatness of functions at a point of $X$ 
(see \ref{flat}). 
This is nothing but a sharpened form of Spallek's theorem. 

We consider that our study still leaves 
a major portion of the extension problem open. 
%We consider that our study is a beginning of the problems. 
We treat only rather easy phenomena in the following sense. 
A set satisfying the 
condition in the main theorem \ref{ext} has always 
the full higher order paratangent space at the accumulating point. 
The construction of the paratangent bundle requires repetition of 
Glaeser operations in general (cf. \cite{bmp2}, $1.8$). 
In contrast, the sets treated in \ref{ext} call for it only once. 
By the result \cite{bmp2}, ${1.3}$, it might be inevitable to 
assume the graphic condition of $\tau_N^e(f)$ with $e>d$ 
for $C^d$ extension in general. 

In the case of the flatness problem, it is interesting to 
analyse the growth of Spallek functions 
(cf. \cite{jam}, ${2.9}$, ${3.6}$, ${4.4}$). 
But we have no idea to connect our present method 
to observe them. 

This work was partly done during the stay at 
Universite des Science et Technologies de Lille. 
The author wish to express his sincere thanks for the courtesy and 
for the helpful discussions with the participants of 
S\'{e}minaire d'Analyse Complexe et Diff\'{e}rentielle of Lille, 
in particular to Professor Ann-Marie Chollet. 
He also would like to thank the participants of 
Seminar of Functions of Complex Variables of Kyoto. 
%%%%%%%%%%%%%%%%%%%%%%%%%%%%%%%%%%%%%%%%%%%%%%%%%%%%%
\section{Multivariate Lagrange interpolation}
%%%%%%%%%%%%%%%%%%%%%%%%%%%%%%%%%%%%%%%%%%%%%%%%%%%

Let us recall some elementary facts on Lagrange interpolation 
in $\mathbb{R}^n$. In the following, functions are $\mathbb{R}$ 
valued and linearity is over $\mathbb{R}$. 
%%%%%%%%%%%%%%%%%%%%%
\begin{proposition}\label{first}
Let $A$ be a subset of $\mathbb{R}^n$ of 
$N$ distinct points 
and $f_1$,\ldots,$f_N$ be functions defined on $A$. 
Then the following conditions are equivalent. 
\begin{enumerate}
\item
$A$ is not contained in the vanishing locus of any non-trivial 
linear combination of $f_1$,\ldots,$f_N$. 
\item
For any set of values prescribed at each point of $A$, 
there exists at most one linear combination of $f_1$,\ldots,$f_N$ 
which takes these values at each point of $A$. 
\item
For any set of values prescribed at each point of $A$, 
there exists at least one linear combination of $f_1$,\ldots,$f_N$ 
which takes these values at each point of $A$. 
\end{enumerate}
\end{proposition}
%%%%%%%%%%%%%%%%%%%%%
\textit{Proof.} 
Let $V:=\bigl( f_i(\BFIT{a}_j) \bigr)$ denote 
the square matrix of the values 
of $f_i$ at points $\BFIT{a}_1,\ldots,\BFIT{a}_N$. Then, 
(1) and (2) are equivalent to the condition that the row vectors of 
$V$ are independent. The condition (3) is equivalent to saying that 
the vectors generate the whole $N$ dimensional space. 
Since $V$ is square, (1), (2), (3) are all equivalent to 
the condition that $V$ is regular. 
\QED
%%%%%%%%%%%%%%%%%%
\begin{proposition}\label{choose-finite}
If $S \subset \mathbb{R}^n$ is not contained in the vanishing locus 
of any non-trivial linear combination of $f_1$,\ldots,$f_M$, 
there exists $A \subset S$ such that $\# A = M$ and $A$ is 
not contained in such a locus either. 
\end{proposition}
%%%%%%%%%%%%%%%%%%%%
\textit{Proof.} 
Let 
$W:=\bigl( f_i(\BFIT{s}) \bigr)_{i=1,\ldots,M;\ \BFIT{s}\in S}$ denote 
the (possibly infinite) matrix of the values 
of $f_i$ at points $\BFIT{s} \in S$. 
The rows of $W$ are linearly independent by our assumption. 
Hence, there is an $M \times M$ regular minor matrix. 
Then the set $A$ of the points corresponding to the columns 
of the minor satisfies the condition. 
\QED
%%%%%%%%%%%%%%%%%%%%%%%%%%%%%%%%%%%%%%%%%%%%%%%%%
\section{Polynomial interpolation}\label{polynomial-int}
%%%%%%%%%%%%%%%%%%%%%%%%%%%%%%%%%%%%%%%%%%%%%%%%

For a subset $A$ of $\mathbb{R}^n$, 
let $\Hdeg(A)$ denote the minimum of the degrees of 
non-zero polynomials vanishing on $A$. 
If there is no such polynomial, we put $\Hdeg(A) = \infty$. 
The dimension of the vector space of homogeneous polynomials 
in $n+1$ variables of degree $d$ coincides with 
that of the vector space of polynomials in $n$ variables of 
degree less than or equal to $d$. We express it by 
$$
N(n,d)=N(d,n):=\binom{n+d}{d}=\frac{(n+d)!}{n!d!}.
$$
%%%%%%%%%%%%%%%%%%
\begin{proposition}\label{degree}
If $A \subset \mathbb{R}^n$ and $\Hdeg(A) 
\geq d+1$, then $\# A \geq N(n,d)$. 
\end{proposition}
%%%%%%%%%%%%%%%%%%%%
\textit{Proof.} 
Consider the vectors of dimension $\# A$ whose components are 
the values of a monomial of degree less than or 
equal to $d$ at 
points of $A$. Such vectors are $N(n,d)$ in number. 
If $\# A < N(n,d)$, these vectors are linearly dependent. 
This implies that some non-trivial linear combination of 
the monomials vanishes at each point of $A$. 
This contradicts the assumption $\Hdeg(A) \geq d+1$. 
\QED
Let us take the set of multi-indices 
$$
I=I(n,d)
:=\bigl\{ \BFIT{i}:=(i_1,\ldots,i_n) \in \{ 0,\ldots,d \}^n: 
|\BFIT{i}| \le d \bigr\}
$$$$
(|\BFIT{i}|:=i_1+\cdots +i_n)
$$
and express the monomials and the derivatives in 
$\BFIT{x}:=(x_1,\ldots,x_n)$ as follows: 
$$
\BFIT{x}^{\SCBFIT{i}}:=x_1^{i_1}\cdots x_n^{i_n}
\quad
\bigl(\BFIT{i}:=(i_1,\ldots,i_n) \in I\bigr),
$$$$
f^{(\SCBFIT{p})}
:=\frac{\partial^{|\SCBFIT{p}|}f}%
{\partial^{p_{1}} x_1\cdots \partial^{p_{n}} x_n}
\quad
\bigl(\BFIT{p}:=(p_1,\ldots,p_n)\in I\bigr).
$$
%%%%%%%%%%%%%%%%%%%%%
\begin{lemma}\label{formula}
The sum of the degrees of all monomials 
in $n$ variables of degrees less than or equal to $d$ is 
equal to the following numbers. 
$$
  \sum_{\SCBFIT{p} \in I(n,d)} |\BFIT{p}|
= n \sum_{i=1}^d i \cdot N(n-1,d-i) 
= \sum_{i=1}^d i \cdot N(n-1,i) 
$$$$
= n \sum_{i=0}^{d-1} N(n,i) 
= n \cdot N(n+1,d-1).
$$ 
\end{lemma}
%%%%%%%%%%%%%%%%%%%%%
\textit{Proof.} 
Let $x_1,\ldots,x_n$ be the variables. 
Let $S_k$ $(k=1,2,\ldots,5)$ denote the $k$-th expression 
in the equality above. 
The first expression $S_1$ is just the quantity 
mentioned at the top of the lemma. The expression 
$S_2$ is obtained by counting the degrees in $x_p$ separately 
for each $p$. 
The summand is the sums for the terms of degree just $i$ in $x_p$. 
The expression 
$N(n-1,d-i)$ denotes the number of such terms, of degree 
equal to or smaller than $d-i$ in the variables other than $x_p$. 
The preceding multiplier $n$ is the number of the choice of $p$. 
The summand of $S_3$ is equal to the product of the degree $i$ and 
the number of the monomial bases 
in $n$ variables of degree just $i$. 
The equality $S_3 = S_4$ follows from 
$$
   \sum_{i=1}^d i \cdot N(n-1,i)
=  \sum_{i=1}^d \frac{(n+i-1)!}{(n-1)!(i-1)!}
=  n \sum_{i=0}^{d-1} N(n,i).
$$
The last equality $S_4=S_5$ follows from the obvious equality 
$$
N(n,i)=N(n+1,i)-N(n+1,i-1).
$$ 
\QED
Suppose that $A$ is a set of 
$N(n,d)$ distinct points in $\mathbb{R}^n$ indexed as 
$$
A:=\bigr\{ \BFIT{a}_{\SCBFIT{i}}
:=(a_{i_1} ,\ldots, a_{i_n})
: \BFIT{i}\in I \bigr\}.
$$
Fixing an ordering of $I$, 
we obtain an $N(n,d)\times N(n,d)$ matrix 
$$
V(A):=V(\BFIT{a}_{\SCBFIT{i}}:\BFIT{i}\in I ) 
:= (\BFIT{a}_{\SCBFIT{i}}^{\SCBFIT{j}})
= (a_{i_1}^{j_1} \cdots a_{i_n}^{j_n}),
$$
where $\BFIT{i}\in I$ are multi-suffixes seen as the row indices 
and $\BFIT{j}\in I$ are multi-exponents seen as the column indices. 
$V(A)$ is called the $n$-dimensional 
\textit{Vandermonde matrix} of $A$ (cf. [AS]). 
This has the following properties. 
%%%%%%%%%%%%%%%%%%%%%
\begin{proposition}\label{vandermonde}
Let $A :=\{ \BFIT{a}_{\SCBFIT{i}}: \BFIT{i}\in I \}\subset \mathbb{R}^n$ 
be a subset of $N(n,d)$ distinct points. 
\begin{enumerate}
\item
$\Det V(A)$ is homogeneous of degree $N(n+1,d-1)$ 
with respect to the $p$-th coordinates of 
$\BFIT{a}_{\SCBFIT{i}}$ $(\BFIT{i} \in I)$ for each fixed 
$p = 1,\ldots,n$. 
\item
Let $\varphi: \mathbb{R}^n \longrightarrow \mathbb{R}^n$ be a linear 
transformation expressed by a matrix $P$. Then 
$$
\Det V\bigl(\varphi(A)\bigr) = (\Det P)^{N(n+1,d-1)}\Det V(A). 
$$
\item
$\Det V(A)$ is an invariant of translations i.e.
$$
\Det V(A-\BFIT{v}) = \Det V(A)\quad (\BFIT{v} \in \mathbb{R}^n).
$$
\end{enumerate}
\end{proposition}
%%%%%%%%%%%%%%%%%%%%%%%
\textit{Proof.} 
\begin{enumerate}
\item
This follows from \ref{formula} ($S_1=S_5$ in the proof).
\item
The general linear group $GL(n,\mathbb{R})$ is generated by 
the transformations of the following forms. 
\begin{enumerate}
\item 
$y_1=x_1,\ldots,y_p=x_p+\lambda x_q,\ldots,y_n=x_n\quad(p \neq q)$,
\item 
$y_1=x_1,\ldots,y_p=\lambda x_p,\ldots,y_n=x_n$ $(\lambda\neq 0)$.
\end{enumerate}
In the case of (a), $\Det P=1$ and $\Det V(\varphi(A)) = \Det V(A)$. 
In the case of (b), $\Det P=\lambda$. Since $\Det V(A)$ is 
homogeneous of degree $N(n+1,d-1)$ 
with respect to the $p$-th coordinates of 
$\BFIT{a}_{\SCBFIT{i}}$ $(\BFIT{i} \in I)$, we have 
$$
\Det V(\varphi(A)) 
= \lambda^{N(n+1,d-1)}\Det V(A)
$$$$
= (\Det P)^{N(n+1,d-1)}\Det V(A).
$$
These prove the equalities.
\item
In view of (2), 
we have only to prove the invariance with respect to the 
transformations in $x_1$ direction, 
which follows in the same way as the case (a) above. 
\end{enumerate}
\QED
%%%%%%%%%%%%%%%%%%%%%%%%
\begin{remark}\label{signed}
{\upshape 
We call 
$\varphi: \mathbb{R}^n \longrightarrow \mathbb{R}^n$ 
\textit{affine} if it is a composition of 
an invertible linear transformation $\varphi'$ 
after (or before) a translation $\varphi''$. 
If $\varphi$ is an affine transformation, the determinant of its 
linear part $\varphi'$ is the ratio of the signed volumes of 
$\varphi(B)$ and $B$ for any measurable subset $B$. 
Thus the proposition above implies that the quotient 
$V(\varphi(A))/V(A)$ 
is expressed as the $N(n+1,d-1)$-th power of the ratio of the 
signed volumes of $\varphi(B)$ and $B$. 
}
\end{remark}
%%%%%%%%%%%%%%%%%%%%%
Applying \ref{first} and \ref{degree} 
to the monomials $\BFIT{x}^{\SCBFIT{j}}$, we have the following. 
%%%%%%%%%%%%%%%%%%%%%
\begin{proposition}\label{position}
Suppose that $A$ is 
a set of $N(n,d)$ distinct points in $\mathbb{R}^n$ with 
$n,\ d \in \mathbb{N}$. 
Then $\Hdeg(A)\le d+1$ and the following conditions are equivalent. 
\begin{enumerate}
\item
$\Hdeg(A)=d+1$. 
\item
$\Det V(A)\neq 0$.
\item
There is a unique polynomial of degree at most $d$ 
which takes the set of values prescribed at each point of $A$. 
\end{enumerate}
\end{proposition}
%%%%%%%%%%%%%%%%%%%%%
The author knows the following geometric interpretations 
of these conditions. 
\begin{enumerate}
\item
The case $n=1$ is well known: $\Det V(A)\neq 0$. 
\item
The case $d=1$ implies that $n+1$ points 
$\BFIT{a}_i:=(a_{i1},\ldots,a_{in})$ $(i=0,\ldots,n)$ in $\mathbb{R}^n$ 
are contained in a hyperplane if and only if 
$$
\left| \begin{array}{cccc}
a_{01} & \ldots & a_{0n} & 1\\
\hdotsfor{4}\\
a_{n1} & \ldots & a_{nn} & 1
\end{array} \right|
= 0.
$$
\item
Suppose that $n=d=2$. 
Let $\BFIT{b}_1$, $\BFIT{b}_2$, $\BFIT{b}_3$ 
denote the intersection points of 

the line joining $\BFIT{a}_{110}$ and $\BFIT{a}_{020}$ 
and 
the line joining $\BFIT{a}_{101}$ and $\BFIT{a}_{002}$; 

the line joining $\BFIT{a}_{011}$ and $\BFIT{a}_{002}$ 
and 
the line joining $\BFIT{a}_{110}$ and $\BFIT{a}_{200}$; 

the line joining $\BFIT{a}_{101}$ and $\BFIT{a}_{200}$ 
and 
the line joining $\BFIT{a}_{011}$ and $\BFIT{a}_{020}$ 

\noindent respectively. By Pascal's theorem, 
(1) in the proposition is equivalent to the condition that 
$\BFIT{b}_1$, $\BFIT{b}_2$, $\BFIT{b}_3$ are not contained in a line. 
\item
A sufficient condition for (2) for general $n$, $d$ 
is given in \cite{as}, ${3.6}$. 
A convenient version of this condition is given in \cite{jam}. 
(The author does not know whether this convenient 
version loses generality or not in comparison to 
\cite{as}, ${3.6}$.) These assure that a general set $A$ 
of $N(n,d)$ points has non-vanishing $\Det V(A)$. 
\end{enumerate}
%%%%%%%%%%%%%%%%%%%%%%%%%%%%%%%%%%%%%%%%%%%%%%%%%%%%%
\section{Higher order paratangent bundle}\label{ptg}
%%%%%%%%%%%%%%%%%%%%%%%%%%%%%%%%%%%%%%%%%%%%%%%%%%%%

Bierstone, Milman and Paw\l ucki \cite{bmp2} have defined the 
\textit{higher order paratangent bundle} $\tau_N^d(X)$ 
of order $d$ for a subset $X$ of a Euclidean space 
(or of a manifold), generalising Glaeser's paratangent bundle. 
We briefly describe the necessary part here. 
Note that we adopt the general definition of $\tau_{N,\SCBFIT{x}}^d(X)$ 
explained in the last section of \cite{bmp2}, whereas 
they used only the case $N=1$ in the 
main part of the paper. Necessity of using larger $N$ was 
already pointed out and related lemmas were prepared by them. 

Let $X$ be a metric space and $V$ a finite dimensional $\mathbb{R}$ 
vector space. 
We call a subset $E \subset X\times V$ a \textit{bundle} 
(of subspaces of $V$) over $X$ 
if the \textit{fibres} $E_{\SCBFIT{a}}:=\{ v: (\BFIT{a},v)\in E \}$ 
are linear subspaces of $V$. 

Let $X$ be a subset of a Euclidean space $\mathbb{R}^n$. 
Let $\CAL{P}_d$ be the $\mathbb{R}$ vector space of polynomials 
of degree equal to or less than $d$ 
and $\CAL{P}^*_d$ its dual vector space. 
Let us denote assignment of the value 
$(-1)^{|\SCBFIT{p}|}
({\partial^{|\SCBFIT{p}|}f}/{\partial\BFIT{x}^{\SCBFIT{p}}})(\BFIT{a})$
to $f \in \CAL{P}_d$ by $\delta^{(\SCBFIT{p})}_{\SCBFIT{a}}$. 
This expresses the derivative of the Dirac delta function 
of order $\BFIT{p}$ with $|\BFIT{p}| \le d$. 
We adopt this symbol because of notational simplicity. 
But take care that it has different properties 
according to $d$. For example, $\delta_a=\delta_b$ always holds 
in $\CAL{P}^*_0$ but not in $\CAL{P}^*_1$. We have 
$\delta_a=\delta_b-\sum(a_i-b_i)\delta^{(e_i)}_b$ and 
$\delta^{(e_i)}_a=\delta^{(e_i)}_b$ in $\CAL{P}^*_1$, 
where $e_i$ is the multi-index whose $i$-th component is 1 and 
other components are 0. In general, the derivatives 
${(-1)^{|\SCBFIT{p}|}}\delta_{\SCBFIT{a}}^{(\SCBFIT{p})}/{\BFIT{p}!}$ 
$(\BFIT{p} := p_1!,\ldots,p_n!,\ |\BFIT{p}| \le d)$ 
form the dual basis of $\{ (\BFIT{x}-\BFIT{a})^{\SCBFIT{p}}:\ 
|\BFIT{p}| \le d \}$ 
and all the derivatives of $\delta_{\SCBFIT{a}}$ 
are expressed by those of $\delta_{\SCBFIT{b}}$ in $\CAL{P}^*_d$ 
(see the first equality in the proof of \ref{ext}). 

Let us put 
$$
E_0 := \bigl\{ (\BFIT{a},\lambda\delta_{\SCBFIT{a}}) :\ 
\BFIT{a}\in X,\ \lambda\in\mathbb{R} \bigr\}.
$$
We define $E_k$ inductively as follows. If $E_k$ is defined, put 
$$
\Delta E_k 
:= \Bigl\{ (\BFIT{a}_0,\ldots,\BFIT{a}_N,\xi_0+\cdots+\xi_N) :\ 
\BFIT{a}_i\in X,\ \xi_i \in E_{k,\SCBFIT{a}_i},
$$$$
|{\BFIT{a}}_i-{\BFIT{a}}_0|^{d-|\alpha|}
\bigm|\xi_i
\bigl(
({\BFIT{x}}-{\BFIT{a}}_i)^{\alpha}
\bigr)\bigm| 
\le 1\ (|\alpha|\le d,\ 0\le i\le N) \Bigr\}
$$
and 
$$
E'_k := \pi\Bigl(\overline{\Delta E_k} 
\cap \bigl\{ (\BFIT{a},\ldots,\BFIT{a},\xi) :\ \BFIT{a}\in X,\ 
\xi\in\CAL{P}^*_d \bigr\}\Bigr),
$$
where 
$\pi:\ X\times\cdots\times X\times \CAL{P}^*_d\longrightarrow X
\times \CAL{P}^*_d$ 
denotes the canonical projection onto the first $X$ 
times $\CAL{P}^*_d$. The intersection of 
the closure of $\Delta E_k$ and 
the diagonal 
coincides with the set of all the limiting points of 
$$
(\BFIT{a}_0,\ldots,\BFIT{a}_N,\xi_0+\cdots+\xi_N) \in \Delta E_k
$$ 
when $\BFIT{a}_0,\ldots,\BFIT{a}_N$ approach to $\BFIT{a}$. 
Finally we put 
$$
E_{k+1} 
:= \bigcup_{\SCBFIT{a}\in X} \bigl(\{ \BFIT{a} \}\times 
\Span E'_{k,\SCBFIT{a}}\bigr)
\subset X\times \CAL{P}^*_d,
$$
where $\Span E'_k$ denotes the linear span of $E'_k$ in the 
fibre. 
The procedure of obtaining $E_{k+1}$ from $E_k$ is an example of 
\text{Glaeser operation} in \cite{bmp2}. 
The sequence $E_1 \subset E_2 \subset E_3 \subset \cdots$ 
stabilizes and we have 
$E_k=E_{2\dim \CAL{P}^*_d}$ $(k \geq 2\dim \CAL{P}^*_d)$ 
as a general property of Glaeser operation 
(\cite{g1}; \cite{bmp2}, ${3.3}$). 
This saturation $\tau_N^d(X) := E_{2\dim \CAL{P}^*_d}$ is called 
the \textit{paratangent bundle} of order $d$ of $X$ 
(\cite{bmp2}). This is a closed subbundle 
of $X\times \CAL{P}^*_d$ in the obvious sense. 
It is known that, for a subbundle (of subspaces), 
closedness is equivalent to upper semi-continuity of inclusion 
(\cite{choquet}, p.67).
Let us call the fibre $\tau_{N,\SCBFIT{a}}^d(X)$ 
the \textit{paratangent space} of order $d$ of $X$ at $\BFIT{a}$. 
Glaeser's (linearized) paratangent bundle is isomorphic 
to $\tau_1^1(X)$. 
%%%%%%%%%%%%%%%%%%%%%%%%%%%%%%%%
\begin{remark}
{\upshape
We can replace the control condition 
$$
|\BFIT{a}_i-\BFIT{a}_0|^{d-|\alpha|}\bigm|\xi_i
\bigl(
(\BFIT{x}-\BFIT{a}_i)^{\alpha}
\bigr)\bigm| 
\le 1\quad(N\le r)
$$
by 
$$
|\BFIT{a}_i-\BFIT{a}_0|^{d-|\alpha|}\bigm|\xi_i
\bigl((\BFIT{x}-\BFIT{a}_i)^{\alpha}\bigr)\bigm| 
\le c
$$
or by
$$
|\BFIT{a}_i-\BFIT{a}_0|^{d-|\alpha|}\bigm|\xi_i
\bigl((\BFIT{x}-\BFIT{a}_0)^{\alpha}\bigr)\bigm| 
\le c
$$
with any $c>0$ independent of $i$ and $\alpha$ (\cite{bmp2}, \S 5). 
This control condition is used to prove the easy half of the 
conjecture below 
(cf. $(\BF{4.7})$, $(\BF{4.16})$, $(\BF{4.17})$ of \cite{bmp2}). 
Hence the control condition seems to endow 
$\tau_{N,\SCBFIT{a}}^d(X)$ a character peculiar to class $C^d$. 
}
\end{remark}
%%%%%%%%%%%%%%%%%%%%%%%%
Now we describe the construction of $\nabla^d_N f$ 
in order to have its image, although we do not use its 
explicit form later. 
Consider a continuous function $f:\ X\longrightarrow \mathbb{R}$ and 
the bundle 
$$
\Phi_0 := \bigl\{ (\BFIT{a},\lambda\delta_{\SCBFIT{a}},\lambda f(\BFIT{a})):\ 
\BFIT{a}\in X,\ \lambda\in\mathbb{R} \bigr\}
\subset X \times \CAL{P}^*_d \times \mathbb{R}
$$
over $X$. If $\Phi_k$ is defined, put 
$$
\Delta \Phi_k := \Bigl\{ (\BFIT{a}_0,\ldots,\BFIT{a}_N,\xi_0+\cdots+\xi_N,
\lambda_0+\cdots+\lambda_N) :
$$$$
\BFIT{a}_i\in X,\ 
(\xi_i,\lambda_i) \in \Phi_{k,a_i},\ 
|\BFIT{a}_i-\BFIT{a}_0|^{d-|\alpha|}\bigm|\xi_i
\bigl(
(\BFIT{x}-\BFIT{a}_i)^{\alpha}
\bigr)\bigm| \le 1
$$$$
(|\alpha|\le d,\ 0\le i\le N) \Bigr\}
$$
and 
$$
\Phi'_k := \pi\Bigl(\overline{\Delta \Phi_k} 
\cap \bigl\{ (\BFIT{a},\ldots,\BFIT{a},\xi,\lambda) :\ 
\BFIT{a}\in X,\ \xi\in\CAL{P}^*_d,\ \lambda\in\BB{R} \bigr\}\Bigr),
$$
where 
$$
\pi:\ X \times\cdots\times X\times \CAL{P}^*_d\times \mathbb{R}
\longrightarrow X\times \CAL{P}^*_d\times \mathbb{R}
$$
denotes the canonical projection onto the first $X$ 
times $\CAL{P}^*_d\times \BB{R}$. Finally we put 
$$
\Phi_{k+1} 
:= \bigcup_{\SCBFIT{a}\in X} \bigl(\{ \BFIT{a} \}\times 
\Span \Phi'_{k,\SCBFIT{a}}\bigr)
\subset X\times \CAL{P}^*_d\times\mathbb{R}.
$$
Since the extension $\Phi_k\subset\Phi_{k+1}$ 
is also a Glaeser operation, 
the sequence $\Phi_0\subset\Phi_1\subset\Phi_2\subset\ldots$ 
stabilizes for $i \geq 2\dim (\CAL{P}^*_d\times \mathbb{R})$
(\cite{g1}, \cite{bmp2}, ${3.3})$ and 
the saturation is denoted by $\nabla^d_N f$. 
%%%%%%%%%%%%%%%%%%%%%%%
\begin{conjecture}[Bierstone-Milman-Paw\l ucki]\label{conjecture} 
\upshape{
Let $X$ be a closed subset of $\mathbb{R}^n$. 
Then there exists $N \in \BB{N}$ such that 
a continuous function $f:\ X\longrightarrow \mathbb{R}$ 
can be extended to a $C^d$ function if and only if 
$\nabla^d_N f \subset \tau_N^d(X)\times \CAL{P}^*_d$ is 
a graph of a map of $\tau_N^d(X)$ into $\mathbb{R}$.
}
\end{conjecture}
%%%%%%%%%%%%%%%%%%%%%%
Only-if part was proved by themselves (\cite{bmp2}, ${4.17})$. 
They proved if part when $X$ is the closure of 
an open subset of a regular submanifold. The paratangent bundle 
$\nabla_N^d f$ with $N=1$ is sufficient for their proof. 
We prove a sharper form of this result in the next section. 
They also proved the case of compact 
subanalytic sets with some loss of 
differentiability using their deep results \cite{bmp1} 
on composite functions. 
%%%%%%%%%%%%%%%%%%%%%
\begin{remark}\label{everywhere}
\upshape{
Suppose that $X$ is a closed subset of $\mathbb{R}^n$ 
and $Y$ a dense subset of $X$. 
If $\tau_{N,\BF{\xi}}^d(X) = \CAL{P}_d^*$ for any $\BF{\xi} \in Y$, 
then $\tau_{N,\BF{\xi}}^d(X) = \CAL{P}_d^*$ for any 
$\BF{\xi} \in X$ by the closedness of $\tau_{N,\BF{\xi}}^d(X)$. 
Then the conjecture \ref{conjecture} is affirmative 
for such an $X$ by \cite{bmp2}, Proof of ${4.20}$.
}
\end{remark}
%%%%%%%%%%%%%%%%%%%%%%
\begin{remark}\label{mfd}
\upshape{
(\cite{bmp2}, ${4.23}$) 
Suppose that $X \subset M \subset \mathbb{R}^n$, 
where $M$ an $m$-dimensional regular submanifold of 
$\mathbb{R}^n$ and that 
$\tau_{N,\BF{\xi}}^d(X) = \tau_{N,\BF{\xi}}^d(M)$ 
for any $\BF{\xi} \in X$. Then the conjecture 
\ref{conjecture} is affirmative for such an $X$.
}
\end{remark}
%%%%%%%%%%%%%%%%%%%%%%%%%%%%%%%%%%%%%%%%%%%%%%%%%%%%%
\section{Set germs with full higher order paratangent spaces}
%%%%%%%%%%%%%%%%%%%%%%%%%%%%%%%%%%%%%%%%%%%%%%%%%%%%
\begin{theorem}\label{ext}
Let $\{ r_k \} \subset \mathbb{R}$ be a positive sequence and 
$A_k := \{ \BFIT{a}^k_0,\ldots,$
$\BFIT{a}^k_N \} \subset \mathbb{R}^n$ 
sets of $N(n,d)=N+1$ distinct points $(k \in \mathbb{N})$. 
Suppose that 
\begin{enumerate}
\item
$A_k$ is contained in the closed ball of radius $r_k$ centred 
at $\BFIT{a}^k_0$$;$
\item
$\lim_{k\to\infty} \BFIT{a}^k_0 = \BFIT{\xi}$$;$
\item
$\lim_{k\to\infty}r_k=0$$;$
\item
there exists $c>0$ independent of $k$ such that 
$$
|\Det V(A_k)| \geq c \cdot r_k^{n\cdot N(n+1,d-1)}
\quad(k \in \mathbb{N}).
$$ 
\end{enumerate}
If $X$ is a closed subset including $\bigcup A_k$, 
then $\tau_{N,\SCBFIT{\xi}}^d(X) = \CAL{P}_d^*$. 
\end{theorem}
%%%%%%%%%%%%%%%%%%%%%
%\vspace{1ex}
\begin{center}
\unitlength 1mm
\begin{picture}(100,47)(1,1)
\thicklines
\put(11,21){$\Bul$}
\put(9.8,23){$\boldsymbol{\xi}$}
\path(77.5,24)(91,39)
\put(86,30){$r_k$}
\put(77.5,24){$\Bul$}
\put(76,19){$\BFIT{a}_0^k$}
\put(67,23){$\Bul$}
\put(66,25){$\BFIT{a}_1^k$}
\put(79,33){$\Bul$}
\put(76,35){$\BFIT{a}_2^k$}
\put(91,17){$\Bul$}
%\put(90,16){$\BFIT{a}_2^k$}
%%
\put(83,7){$\Bul$}
\put(83,10){$\BFIT{a}_N^k$}
\put(78,24){\circle{40}}
\put(75,8.5){$A_k$}
%%%%%%%%%%%%%%%%%%%%%%%
\put(39,15){\circle{20}}
\put(39,15){$\Bul$}
\put(31,14){$\Bul$}
\put(37,23){$\Bul$}
\put(45,11){$\Bul$}
\put(41,13){$\Bul$}
\put(37.5,17){$\BFIT{a}_0^{k+1}$}
\put(32,8.5){$A_{k+1}$}
%%%%%%%%%%%%%%%%%%%%%%%%
\put(25,27){\circle{10}}
\put(22,29){$\Bul$}
\put(24,25){$\Bul$}
\put(22,26){$\Bul$}
\put(29,26){$\Bul$}
\put(25,27){$\Bul$}
%%%%%%%%%%%%%%%%%%%%%%%%
\put(21,22){\circle{6}}
\put(20,24){$\Bul$}
\put(21,21.7){$\Bul$}
\put(22,20){$\Bul$}
\put(23,23){$\Bul$}
\put(19,22){$\Bul$}
%%%%%%%%%%%%%%%%%%%%%%%%
\put(14,19){\circle{3}}
\put(14,19){$\Bul$}
\put(14.8,18.2){$\Bul$}
\put(15,19.5){$\Bul$}
\put(13.5,20){$\Bul$}
\put(13,18.3){$\Bul$}
%%%%%%%%%%%%%%%%%%%%%%%%
\end{picture}
\end{center}
%\vspace{1ex}
%%%%%%%%%%%%%%%%%%%%%%
\begin{remark}
\upshape{
In view of \ref{position}, 
the condition (4) implies that the points of $A_k$ are algebraically 
in general position in the balls of (1) uniformly with respect to 
$k$. If the interior of $X$ is adherent to 
0, this condition is satisfied. Hence, 
by \cite{bmp2}, Proof of $\BF{4.20}$ 
(or \ref{everywhere} with $Y=X$), we see that 
our theorem is an improvement of \cite{bmp2}, ${4.19}$: 
balls are replaced by the sets $A_k$ of $N+1$ points.
}
\end{remark}
%%%%%%%%%%%%%%%%%%%%%

\textit{Proof of 4.1.} 
Let $\BFIT{p}_0,\ldots,\BFIT{p}_N$ denote the elements of 
$I(n,d)$. If $f$ is a polynomial of degree $d$, we have 
$$
V(A_k-\BFIT{a}^k_{\BF{0}})\,
\left( \begin{array}{c}
\frac{f^{(\SCBFIT{p}_0)}(\BFIT{a}^k_0)}{\BFIT{p}_0!}\\
\vdots\\
\frac{f^{(\SCBFIT{p}_N)}    (\BFIT{a}^k_0)}{    \BFIT{p}_N!}
\end{array} \right)
=
\left( \begin{array}{c}
f(\BFIT{a}^k_0)\\
\vdots\\
f(\BFIT{a}^k_N)
\end{array} \right).
$$
Namely, the Dirac deltas 
$\delta_{\SCBFIT{a}^k_0},\ldots,\delta_{\SCBFIT{a}^k_N}$ are 
expressed in terms of the higher order derivatives 
$\delta_{\SCBFIT{a}^k_0}^{(\SCBFIT{p}_0)},\ldots,
\delta_{\SCBFIT{a}^k_0}^{(\SCBFIT{p}_N)}$ of the Dirac delta at $\BFIT{a}_0^k$ 
in $\CAL{P}_d^*$. 

Since the Vandermonde determinant does not vanish, 
$\delta^{(\SCBFIT{p}_{\SCBFIT{i}})}_{\SCBFIT{a}_0^k}$ are spanned 
by $\delta_{\SCBFIT{a}^k_1}, \ldots , \delta_{\SCBFIT{a}^k_N}$. 
All the elements of 
$V(A_k-\BFIT{a}^k_{\BF{0}})$ with column index $i$ 
are homogeneous polynomials of degree $|\BFIT{p}_i|$  
in all the components of all $\BFIT{a}_i-\BFIT{a}_0$ and the 
Vandermonde determinant $\Det V(A_k-\BFIT{a}^k_{\BF{0}})$ is 
a homogeneous polynomial of degree $n\cdot N(n+1,d-1)$ 
in them. Then the elements with row index $i$ of 
the cofactor matrix 
$\bigl(\Det V(A_k-\BFIT{a}^k_{\BF{0}})\bigr)
V(A_k-\BFIT{a}^k_{\BF{0}})^{-1}$ 
are homogeneous of degree 
$$
\sum_{j\neq i}|\BFIT{p}_j|
=\sum_j|\BFIT{p}_j|-|\BFIT{p}_i|
=n\cdot N(n+1,d-1)-|\BFIT{p}_i|
$$
in them. 
Here, the last equality follows from \ref{formula}. 
Applying the condition (4) and \ref{formula}, we see that 
all the elements of $V(A_k-\BFIT{a}^k_{\BF{0}})^{-1}$ 
with row index $i$ are majorized by a constant 
multiple of $r_k^{-|\SCBFIT{p}_i|} \le r_k^{-d}$. 
Hence the control conditions for the coefficients of 
$\delta_{\SCBFIT{a}^k_0}$ 
in the construction of $\tau_{N,\BF{0}}^d(X)$ 
are satisfied. Then 
$\delta^{(\SCBFIT{p}_j)}_{\BF{0}} \in \tau_{N,\BF{0}}^d(X)$ 
follows as their limits. 
\QED
%%%%%%%%%%%%%%%%%%%%%%%%%%%%%%%%%%%%%%%%%%%%%%%%%%%%%
\section{Paratangent bundles of self-similar sets}
%%%%%%%%%%%%%%%%%%%%%%%%%%%%%%%%%%%%%%%%%%%%%%%%%%%%

First we recall the definition of self-similar set. 
The readers can refer to \cite{falconer} and \cite{yhk} 
for further explanation. 
A map $\varphi: \mathbb{R}^n \longrightarrow \mathbb{R}^n$ 
is called a 
\textit{contraction} if there exists $K \in (0,1)$ such that 
$$
\|\varphi(\BFIT{x})-\varphi(\BFIT{y})\| 
\le K \cdot \|\BFIT{x}-\BFIT{y}\|
\quad(\BFIT{x},\BFIT{y} \in \mathbb{R}^n).
$$
If a finite set 
$\varphi_1,\ldots,\varphi_p: \mathbb{R}^n 
\longrightarrow \mathbb{R}^n$ 
of contractions of $\mathbb{R}^n$ is given, there exists 
a unique non-empty compact set $S \subset \mathbb{R}^n$ such that 
$S = \bigcup_{i=1}^p \varphi_i(S)$. 
Such an $S$ is called the 
\textit{attractor} or the \textit{invariant set} of 
$\varphi_1,\ldots,\varphi_p$. In particular, every 
contraction $\varphi$ of $\mathbb{R}^n$ has a unique fixed 
point $F(\varphi)$. 
The following is known as \textit{Williams' formula}. 
For the attractor of $\varphi_1,\ldots,\varphi_p$, we have 
$$
S = \overline{\bigcup F(\varphi_{i_1}\COMP\cdots\COMP\varphi_{i_q})}
$$
where $(i_1,\ldots,i_q)$ runs over all finite sequences 
of elements of 
$\{ 1,\ldots, p \}$. Let us call an affine transformation 
a \textit{similarity transformation} 
if it preserves the angle of every ordered triplet of points. 
If $\varphi$ is a similarity transformation, there exists 
$\lambda \in (0,\infty)$ such that 
$$
\|\varphi(\BFIT{x})-\varphi(\BFIT{y})\| 
= \lambda \cdot \|\BFIT{x}-\BFIT{y}\|
\quad(\BFIT{x},\BFIT{y} \in \mathbb{R}^n).
$$
We call $\lambda$ the \textit{similarity ratio} of $\varphi$. 
An attractor of $\varphi_1,\ldots,\varphi_p$ $(p \geq 2)$ is called 
\textit{self-similar}, if 
all $\varphi_i$ are similarity transformations (see \cite{falconer}). 
(Often more general attractors are called self-similar 
(see \cite{yhk}, p.18).)
The next lemma is almost immediate from \ref{signed}.
%%%%%%%%%%%%%%%%%%%%%%
\begin{lemma}\label{ratio}
Let $A$ be a point set with $\# A = N+1 = N(n,d)$ and $B$ 
another point set similar to $A$. Then we have 
$$
\frac{\Det V(A)}{\delta(A)^{n\cdot N(n+1,d-1)}}
=
\frac{\Det V(B)}{\delta(B)^{n\cdot N(n+1,d-1)}},
$$
where $\delta$ denotes the diameter. 
\end{lemma}
%%%%%%%%%%%%%%%%%%%%%
\begin{theorem}\label{self-similar}
For any $d \in \mathbb{N}$, we define $N$ by $N(n,d) = N+1$. 
Let $X \subset \mathbb{R}^n$ be a self-similar subset 
with $\Hdeg(X) \geq d+1$. 
Then $\tau_{N,\SCBFIT{s}}^d(X) = \CAL{P}_d^*$ for any $\BFIT{s} \in X$. 
\end{theorem}
%%%%%%%%%%%%%%%%%%%%%
\newpage% PRF

\textit{Proof.} 
Suppose that $X$ is defined by contracting similarity 
transformations $\varphi_1,\ldots,\varphi_p$. 
By \ref{everywhere}, 
we have only to prove that $\tau_{N,\SCBFIT{s}}^d(X) = \CAL{P}_d^*$ 
for points $\BFIT{s}$ of a dense subset of $X$. Then, by 
Williams' formula, we may assume that $\BFIT{s}$ is the fixed point 
of $\varphi_{i_1}\COMP\cdots\COMP\varphi_{i_k}$, i.e. 
$\BFIT{s} = F(\varphi_{i_1}\COMP\cdots\COMP\varphi_{i_k})$. 

Since $\Hdeg X \geq d+1$, 
there exists an $N+1$ point subset $A \subset X$ 
such that $\Hdeg A \geq d+1$, by \ref{choose-finite}. 
This implies that $V(A) \neq 0$ 
(and $\Hdeg A=d+1)$ by \ref{position}. 
Since $\varphi_{i_1},\ldots,\varphi_{i_k}$ are 
similarity transformations, 
so is $\varphi_{i_1}\COMP\cdots\COMP\varphi_{i_k}$. 
Let $\lambda>0$ denote its 
similarity ratio. Of course, $\lambda < 1$. 

The set $A_k := (\varphi_{i_1}\COMP\cdots\COMP\varphi_{i_k})^k(A)$ 
is included in the closed ball of radius 
$$
\lambda^k \cdot \max\{ |\BFIT{x}-\BFIT{s}|:\ \BFIT{x}\in A \}
$$ 
with centre $\BFIT{s}$. If we number the points of $A_k$ as 
$A_k = \{ \BFIT{a}_0^k,\ldots,\BFIT{a}_N^k \}$ 
arbitrarily, $A_k$ is contained in the ball of radius 
$$
r_k := 2 \lambda^k \cdot \max\{ |\BFIT{x}-\BFIT{s}|:\ \BFIT{x}\in A \} 
$$
centred at $\BFIT{a}_0^k$. We know that 
$$
\frac{\Det V(A_k)}{\delta(A_k)^{n \cdot N(n+1,d-1)}}
= \frac{\Det V(A_1)}{\delta(A_1)^{n \cdot N(n+1,d-1)}}
$$
by \ref{ratio}. Since this expression 
and the ratio $\delta(A_k)/r_k$ are independent of $k$, 
the condition (4) of the theorem holds. 
\QED
Most of the classical fractal sets constructed geometrically 
are self-similar and not contained in an 
algebraic hypersurface (hence $\Hdeg(X) \geq d+1$). 
Among them are Cantor set, 
Koch curve, Sierpinski gasket and Menger sponge. 
Non-algebraicity of these fractal sets follows from 
the fact that the local Hausdorff dimension of a proper 
algebraic subset is smaller than that of the ambient space. 
%%%%%%%%%%%%%%%%%%%%%%%%%%%%%%%%%%%%%%%%%%%%%%%%%%%%%
\setcounter{section}{5}% PRF
\section{Control of flatness by values}
%%%%%%%%%%%%%%%%%%%%%%%%%%%%%%%%%%%%%%%%%%%%%%%%%%%%

Let $f$ be a $C^d$ function defined on an open neighbourhood 
of $0 \in \mathbb{R}^n$. Let us call $f$ $k$\textit{-flat} if 
$f^{(\SCBFIT{p})}(0)=0$ for $\BFIT{p}:=(p_0,\ldots,p_n)\in I(n,d)$ with 
$|\BFIT{p}| \le k \le d$. 
%%%%%%%%%%%%%%%%%%%%%
\begin{remark}
\upshape{
As to the terms \textit{flatness} and \textit{order}, 
the author now understood that it is better to use both 
depending on the category of functions. When we treat 
analytic functions, order is convenient because it is 
a valuation (or related to valuations, on a singular space), 
a familiar notion to algebraists. 
If the order of $f$ is $p$, then $f$ is of course $(p-1)$-flat. 
When we treat $C^d$ functions for finite $d$, 
there occurs a difficulty in defining order. 
If all the partial derivatives of $f$ vanishes order up to $d$, 
$f$ is $d$-flat. But we can not define its order confidently, 
so long as we permit non-integer values. 
So flatness is better in this category. 
}
\end{remark}
%%%%%%%%%%%%%%%%%%%%%%%%%%
\begin{theorem}\label{flat}
Let $p>0$ be a positive number, $\{ r_k \} \subset \mathbb{R}$ 
a positive sequence 
and $A_k := \{ \BFIT{a}^k_0,\ldots,\BFIT{a}^k_N \} \subset \mathbb{R}^n$ 
$(k \in \mathbb{N})$ sets of $N(n,d)=N+1$ distinct points. 
Suppose that: 
\begin{enumerate}
\item
$A_k$ is contained in the closed ball of radius $r_k$ centred 
at $\BFIT{a}^k_0$$;$
\item
$\lim_{k\to\infty}r_k=0$$;$
\item
$\lim_{k\to\infty}\BFIT{a}_0^k=0$$;$
\item
there exist $c,\ e>0$ such that $|\Det V(A_k)| \geq c \cdot r_k^{e}$ 
$(k \in \mathbb{N})$. 
\end{enumerate}
For a $C^{d}$ $(d\ge p)$ 
function $f$ defined in a neighbourhood of $\BF{0}$, 
we put  
$$
S_k := 
r_k^{-p} \cdot 
\max\bigl\{ |f(\BFIT{x})|: \BFIT{x}\in A_k 
\bigr\},
\quad
m := p-\bigl(e-n\cdot N(n+1,d-1)\bigr).
$$
\begin{description}
\item[(i)] 
If $m$ is an integer and $\lim_{k\to\infty}S_k=0$, 
then $f$ is $m$-flat at $\BF{0}$. 
\item[(ii)]
If $m$ is not an integer and $S_k$ is bounded, 
then $f$ is $[m]$-flat at $\BF{0}$, 
where $[m]$ denotes the maximal integer not greater than $m$. 
\end{description}
\end{theorem}
%%%%%%%%%%%%%%%%%%%%%%%%%
\begin{remark}\label{factor}
\upshape{
This theorem is useful in the following situation. 
Let $\{ s_k \} \subset \mathbb{R}$ be a positive sequence. 
Suppose that $A_k$ is contained in the closed ball of radius 
$s_k$ centred at $\BF{0}$. If $p,\ q>0$, $S_k$ in the theorem 
is majorized by 
$$
T_k := 
\frac{s_k^q\ }{r_k^p} \cdot 
\max\Bigl\{ \Bigm|\frac{f(\BFIT{x})}{\BFIT{x}^q}\Bigm|: \BFIT{x}\in A_k 
\Bigr\}.
$$
The first factor of $T_k$ is concerned with 
the shrinking of balls containing $A_k$ 
and the second with flatness of the values of 
$f$ along $\bigcup A_k$. 
If the $T_k$ tend to 0, then so do the $S_k$. 
If the $T_k$ are bounded, then so are the $S_k$. 
}
\end{remark}
%%%%%%%%%%%%%%%%%%%%%%%%%%
\begin{remark}
\upshape{
The expression $m$ above is rather complicated. 
We can understand this as follows. 
If the conditions in the theorem 
holds, then $e$ must satisfy $e \geq n\cdot N(n+1,d-1)$ 
by \ref{formula}. The equality here means that the 
points of each $A_k$ are algebraically in general position 
in the balls of (1) ``uniformly with respect to $k$". 
If this is the case, we have $m = p$ and, in view of 
\ref{factor}, (i) is a sharpening of Spallek's theorem 
\cite{spallek}, {1.4}: balls are replaced by sets $A_k$ 
of $N+1$ points. 
The term $e-n\cdot N(n+1,d-1) \geq 0$ is the adjustment for the 
case when the algebraic genericities of the positions of 
the points of $A_k$ degenerate as $k$ increases. 
}
\end{remark}
%%%%%%%%%%%%%%%%%%%%%%%%%%
\textit{Proof of 6.2.} 
We may assume that $f$ is defined in a neighbourhood of 
the closure of the convex hull of $\bigcup A_k$. 
Let us adopt an ordering of $I(n,d)$ such that 
$$
\{ \BFIT{p}: |\BFIT{p}|<d \} = \{ \BFIT{p}_0,\ldots,\BFIT{p}_M \},
\quad% $$$$
\{ \BFIT{p}: |\BFIT{p}|=d \} = \{ \BFIT{p}_{M+1},\ldots,\BFIT{p}_N \}
$$$$(M=N(n,d-1)).$$
Then by the Taylor formula, 
there exists $\theta_i^k \in (0,1)$ such that we have 
$$
V(A_k-\BFIT{a}^k_{\BF{0}})\,
\left( \begin{array}{c}
\frac{f^{(\SCBFIT{p}_0)}(\BFIT{a}^k_0)}{\BFIT{p}_0!}\\
\vdots\\
\frac{f^{(\SCBFIT{p}_M)}(\BFIT{a}^k_0)}{\BFIT{p}_M!}\\
\frac{f^{(\SCBFIT{p}_{M+1})}(\BFIT{a}_0^k)}{\BFIT{p}_{M+1}!}\\
\vdots\\
\frac{f^{(\SCBFIT{p}_N)}    (\BFIT{a}_0^k)}{    \BFIT{p}_N!}
\end{array} \right)
=
\left( \begin{array}{c}
f(\BFIT{a}^k_0)      +  \delta_1^k\\
\vdots\\
f(\BFIT{a}^k_M)      +  \delta_M^k\\
f(\BFIT{a}^k_{M+1})  +  \delta_{M+1}^k\\
\vdots\\
f(\BFIT{a}^k_N)      +  \delta_N^k
\end{array} \right),
$$
where 
$$
\delta_i^k 
:= \sum_{j=M+1}^N
\bigl(f^{(\SCBFIT{p}_j)}(\BFIT{a}_0^k)-f^{(\SCBFIT{p}_j)}(\BFIT{b}_i^k)\bigr)
(\BFIT{a}_i^k-\BFIT{a}_0^k)^{\SCBFIT{p}_j}/\BFIT{p}_j !,
$$$$
\BFIT{b}_i^k := \theta_i^k\BFIT{a}_0^k+(1-\theta_i^k)\BFIT{a}_i^k,\ 
0 < \theta_i^k < 1
\quad (i=0,\ldots,N;\ k=1,2,\ldots)
$$ 
and $V(A_k-\BFIT{a}_{0}^k)$ is the Vandermonde matrix of 
the translation of $A_k$ by $-\BFIT{a}_{0}^k$. 
As we have seen in the proof of \ref{ext}, 
the elements with row index $i$ of 
$$
\bigl(\Det V(A_k-\BFIT{a}^k_{0})\bigr)V(A_k-\BFIT{a}^k_{0})^{-1}
$$ 
are homogeneous of degree 
$$
\sum_{j\neq i}|\BFIT{p}_{j}|
=n\cdot N(n+1,d-1)-|\BFIT{p}_{i}|
$$
in all the components of all $\BFIT{a}_i-\BFIT{a}_0$. 
Hence there exists $C > 0$ such that 
$$
|f^{(\SCBFIT{p}_{i})}(\BFIT{a}_0^k)|
\le
C \cdot r_k^{m - p - |\SCBFIT{p}_{i}|}
\max\{ |f(\BFIT{a}_i^k)| + |\delta_i^k| :\ 1\le i\le N \}.
$$
Sinse $f$ is of class $C^d$ and since 
$\lim_{k\to\infty} \BFIT{a}_0^k = \lim_{k\to\infty} \BFIT{b}_i^k 
= \BF{0}$, we see that 
$$
\lim_{k\to\infty} r_k^{-d}|\delta_i^k| = 0\ (1\le i\le N).
$$ 
Hence, if $|\BFIT{p}_{i}| \le m$, we have 
$$
\lim_{k\to\infty}r_k^{m - p - |\SCBFIT{p}_{i}|} |\delta_i^k|
 = 0\quad(1\le i\le N).
$$
We have assumed that $|f(\BFIT{a}_i^k)|\le r_k^p S_k$. 

In the case of (i), this implies that 
$\lim_{k\to\infty}|f^{(\SCBFIT{p}_i)}(\BFIT{a}_i^k)|=0$ for $\BFIT{p}_{i}$ 
with $|\BFIT{p}_{i}| \le m.$ Since 
$\lim_{k\to\infty} \BFIT{a}_0^k = \BF{0}$, 
we have $f^{(\SCBFIT{p}_i)}(\BF{0})=0$ for such $\BFIT{p}_{i}$,  
which completes the proof of (i). 

In the case of (ii). let us define $\tilde{p}$ and $\tilde{S}_k$ by 
$$
[m] = \tilde{p}-\bigl(e-n\cdot N(n+1,d-1)\bigr),
\quad
\tilde{S}_k := 
r_k^{-\tilde{p}} \cdot 
\max\bigl\{ |f(\BFIT{x})|: \BFIT{x}\in A_k 
\bigr\}.
$$
Since $\tilde{p}<p$, $\lim_{k\to\infty}\tilde{S}_k = 0$ holds 
and (ii) follows from (i). 
\QED
%%%%%%%%%%%%%%%%%%%%%%%%%%%%%%%%%%%%%%%%%%%%%%%%%%

\vspace{1ex}
%%%%%%%%%%%%%%%%%%%%%%%%%%%%%%%%%%%%%%%%%%%%%%%%
\noindent
Department of Mathematics\\
Kinki University\\
Kowakae Higashi-Osaka 577-8502, Japan\\
{\sf e-mail: izumi@math.kindai.ac.jp}
%%%%%%%%%%%%%%%%%%%%%%%%%%%%%%%%%%%%%%%%%%%%%%%%%%%%%
\end{document}